\documentclass[
    aip,instruments
    amsmath,
    amssymb,
    reprint,
    floatfix,
    nofootinbib
]{revtex4-2}
\usepackage{graphicx}
\usepackage{dcolumn}
\usepackage{bm}

\usepackage[utf8]{inputenc}
\usepackage[T1]{fontenc}
\usepackage{mathptmx}
\usepackage{etoolbox}

\usepackage{bm}
\usepackage{graphicx}
\usepackage{amsthm}
\usepackage{notes2bib}
\usepackage{amsmath}
\usepackage{amssymb}
\usepackage{url}
\usepackage{enumerate}
\usepackage{epstopdf}
\usepackage{color}
\usepackage{mathtools}
\usepackage{setspace}
\usepackage{appendix}
\usepackage[normalem]{ulem}
\usepackage[section]{placeins}
\usepackage{lineno}
\usepackage[caption=false]{subfig}
\usepackage{cases}
\newcommand{\allblack}{\color{black}{}}

\newcommand{\allred}{\color{black}{}}

\makeatletter
\def\@email#1#2{%
\endgroup \patchcmd{\titleblock@produce} {\frontmatter@RRAPformat} {\frontmatter@RRAPformat{\produce@RRAP{*#1\href{mailto:#2}{#2}}}\frontmatter@RRAPformat} {}{} }%
\makeatother
\newcommand{\etal}{\textit{et al}.}

\begin{document}
\title{Finite Invariant Sets with Bridging Points in Logistic IFS}

\author{Hibiki Kato}
\affiliation{Graduate School of Business Administration, Hitotsubashi University, Tokyo 186-8601, Japan}
\affiliation{
    Department of Mathematics, University of Maryland, College Park, MD 20742, USA
}
\author{Tamotsu Onozaki}
\affiliation{Faculty of Economics, Rissho University, Tokyo, Japan}

\author{Yoshitaka Saiki*}
\email{yoshi.saiki@r.hit-u.ac.jp}
\affiliation{Graduate School of Business Administration, Hitotsubashi University, Tokyo 186-8601, Japan}

\author{Yasumasa Sugita}
\affiliation{Graduate School of Business Administration, Hitotsubashi University, Tokyo 186-8601, Japan}

\date{\today}

\begin{abstract}
\allred
We investigate iterated function systems (IFS) that randomly alternate between two non-identical one-dimensional maps as simple models of regime-switching dynamical systems. 
Our primary focus is on finite invariant sets exhibiting “toss-and-catch” dynamics, in which trajectories alternate between fixed points and periodic points of the constituent maps. 
Using two representative types of low-dimensional nonlinear systems—a pair of logistic maps and a combination of logistic and tent maps—we derive exact parameter conditions for several toss-and-catch structures. 
The comparison between these systems reveals two distinct mechanisms for finite invariant sets: one mediated by bridging points that connect invariant structures of different maps, and another generated by nontrivial intersections shared by the maps themselves.
Notably, we identify invariant sets containing bridging points that are not periodic points of either constituent map.
These results demonstrate that switching between simple nonlinear maps can generate invariant structures that do not exist in the constituent systems alone, suggesting a more general mechanism for the emergence of invariant structures in random dynamical systems.
\allblack
   \end{abstract}

\maketitle

\noindent Phenomena in which two regimes alternate over time are common in both natural and social systems, ranging from autonomic nervous system regulation to economic cycles.
Regime switching can be modeled by alternation between different dynamic rules.
In this study, we model such alternation using iterated function systems that randomly apply one of the two nonlinear maps. We specifically examine the IFS constructed from fundamental nonlinear maps: the logistic and tent maps.
Although IFS typically exhibits noisy dynamics, we identify a simple structure, ``toss-and-catch,'' defined as a finite invariant set in which trajectories wander among periodic points of the constituent maps. We derive the precise conditions under which such toss-and-catch structures arise.
Our primary contribution is the demonstration that some structures necessarily involve ``bridging points'' within the toss-and-catch that are not periodic points of either constituent map.

\section{\label{sec:introduction} Introduction}

In the real world, it is common to observe phenomena in which two regimes with distinct dynamic rules alternate over time.
A simple example is the natural transition between day and night.
However, some alternations are nonperiodic, such as the switching between the sympathetic and parasympathetic nervous systems in autonomic regulation, changes in the ruling party in a two-party political system, and shifts in corporate behavior between economic expansions and recessions.
These alternating patterns between two states occur frequently, and modeling such behavior provides a realistic approach to understanding complex real-world phenomena.

Various modeling frameworks have been developed to incorporate regime switching. In particular, some explicitly model the dynamics of the regime variable, as in the Markov switching model~\cite{hamilton1989new}, while others implement state-dependent switching, as in hybrid dynamical systems~\cite{bernardo2008}.
Here, we investigate random iterated function systems (IFS), in which the active map changes randomly over time.

Consider a metric space $\mathcal{M}$ and a collection of two maps:
\begin{equation*}
    \mathcal{F}=\{g,h\}, \quad g,h: \mathcal{M} \to \mathcal{M},
\end{equation*}
which we refer to as an IFS.
In this study, we consider a random dynamical system generated by an IFS composed of two distinct maps, randomly selected for each time step \( n \in \mathbb{N} \).
Let random IFS $\{F_n\}$ be defined as:
\begin{equation}
    F_n =
    \begin{cases}
        g & \text{with probability } p, \\
        h & \text{with probability } 1-p,
    \end{cases}
    \label{eq: IFS}
\end{equation}
where $p$ is constant in the interval $(0, 1)$, and the maps are selected independently and identically at each step. Here, the probability represents a simple exogenous model of the regime-switching dynamics. We will refer to this system as the IFS and study its invariant set. We focus on an {\bf invariant set}~\cite{crovisier_ifs_2006, hutchinson_fractals_1981} $ \Lambda $ of the IFS, which satisfies:
\begin{equation}
    \Lambda = g(\Lambda) \cup h(\Lambda).
\end{equation}

The system defined by \eqref{eq: IFS} is an example of a random dynamical system~\cite{arnold1998random}, which can exhibit properties absent in its deterministic counterparts.
Randomness is often introduced via additive~\cite{sano_reduction_2020} or multiplicative noise applied to variables or parameters. In contrast to perturbation-based approaches, the IFS framework can capture more drastic changes in system dynamics, enabling the explicit modeling of switching behavior.

While studies of IFS have traditionally focused on generating fractals using only contracting maps~\cite{barnsley85}, some research has examined IFSs' dynamics with both contracting and expanding structures~\cite{homburg2025}. For example, Pelikan (1984)~\cite{pelikan84} analyzed IFS constructed from
piecewise linear maps, such as \( g(x) = 2x \bmod 1 \) and
\( h(x) = \tfrac{1}{2}x\), and established the existence of
absolutely continuous invariant measures. Gutierrez \etal~ (1993)~\cite{gutierrez1993logistic} investigated logistic maps under dichotomous noise, revealing changes in bifurcation and invariant measure depending on how the noise is applied.

Regarding logistic IFS, Abbasi \etal~ (2018)~\cite{abbasi2018iterated} investigated an IFS composed of logistic maps with different parameter values, focusing on dynamics near fixed points with distinct stability properties: one contracting and the other expanding. They examined intermittent behavior and the ``synchronization'' of trajectories originating from different initial conditions.
These studies demonstrate that the dynamics of an IFS are nontrivial and cannot be understood solely by analyzing the individual maps in isolation.
\allred 
While switching dynamics are often studied in terms of stability, synchronization, intermittency, or statistical properties of trajectories~\cite{hamilton1989new,bernardo2008,pelikan84,homburg2025,abbasi2018iterated}, comparatively less attention has been paid to the possibility that switching itself may generate new invariant structures. 
In particular, invariant sets of a switching system need not be reducible to invariant sets of the constituent maps individually. 
The present work provides explicit low-dimensional examples of such emergent invariant structures.
\allblack

Motivated by this perspective, we analyze an IFS composed of two nonlinear maps $g$ and $h$ from $[0,1]$ to itself.
Specifically, we consider two types of systems: a pair of logistic maps with different parameter values (logistic IFS) and a combination of a logistic map and a tent map (logistic-tent IFS). The dynamics of these systems cannot be trivially derived from those of the constituent maps and, in most cases, yield more complex trajectories.
In this context, we define the limit set obtained after a long-term iteration as the (numerical) attractor of the IFS.

For deterministic maps, invariant sets typically consist of periodic orbits or chaotic attractors; however, in IFSs, invariant sets are not restricted to such cases.
Let $\Lambda_g$ and $\Lambda_h$ denote the union of all periodic points of $g$ and $h$, respectively.
The first case occurs when $\Lambda$ is entirely contained in $\Lambda_g\bigcup\Lambda_h$, that is,
\begin{equation*}
    \Lambda \setminus \left(\Lambda_{g} \bigcup \Lambda_{h} \right) = \emptyset.
\end{equation*}
The second case arises when $\Lambda$ is not contained in $\Lambda_{g} \cup \Lambda_{h}$, that is,
\begin{equation}
    \Lambda \setminus \left(\Lambda_{g} \bigcup \Lambda_{h}\right) \neq \emptyset.
    \label{eq:bridging}
\end{equation}
We focus in particular on cases in which the invariant set $\Lambda$ of an IFS is finite in the rest of the paper. We refer to a finite invariant set $\Lambda$ as a {\bf toss-and-catch} structure. If a toss-and-catch structure consists
of $n$ points, we call it an $n$-point toss-and-catch~\cite{NOTE}.
Throughout this paper, the toss-and-catch structure includes cases in which the fixed points of the two maps coincide. In particular, when the trivial fixed points of the consistent maps coincide, the resulting singleton invariant set is regarded as a trivial one-point toss-and-catch. A nontrivial example is presented in Sec.~\ref{sec:logistic-tent}.
To systematically identify parameter sets exhibiting toss-and-catch structures, we quantify the ``size'' of the IFS attractor. For each parameter pair, we generate a sufficiently long trajectory after discarding transient behavior and estimate the attractor size as the minimum number of $\varepsilon$-balls required to cover the trajectory. Parameter values corresponding to toss-and-catch structures appear as distinctive curves in parameter space, and the attractor size remains relatively small in the neighborhoods of these curves.

\section{Logistic IFS}\label{Sec:logistic-IFS}
In this section, we consider an IFS in which one of two logistic maps is selected at random at each iteration.
In system~\eqref{eq: IFS}, we set \( g = f_\alpha \) and \( h = f_\beta \), where \( f_\gamma(x) = \gamma x(1 - x) \), with \( \gamma \in [0,4]\).
This construction defines a logistic IFS, which typically exhibits noisy or irregular behavior due to random switching between maps.

We define the logistic IFS $F(x)$ with parameters satisfying \( \alpha < \beta \) as follows:
\begin{equation}
    F(x) =
    \begin{cases}
        g(x) = f_{\alpha}(x) =  \alpha x(1-x) & \text{with } p, \\
        h(x) = f_{\beta}(x) = \beta x(1-x) & \text{with } 1-p,
    \end{cases}
\end{equation}
where we set $p=0.5$ in the numerical examples of Secs.~\ref{Sec:logistic-IFS} and \ref{sec:logistic-tent}.
 Without loss of generality, we assume $\alpha < \beta$  since exchanging $\alpha$ and $\beta$ yields identical statistical dynamics when $p=0.5$.

\begin{nolinenumbers}
    \begin{figure}[h]
        \centering
        \subfloat[deterministic logistic map\label{fig:bif_compare:a}]{
            \includegraphics[width=0.696\columnwidth]{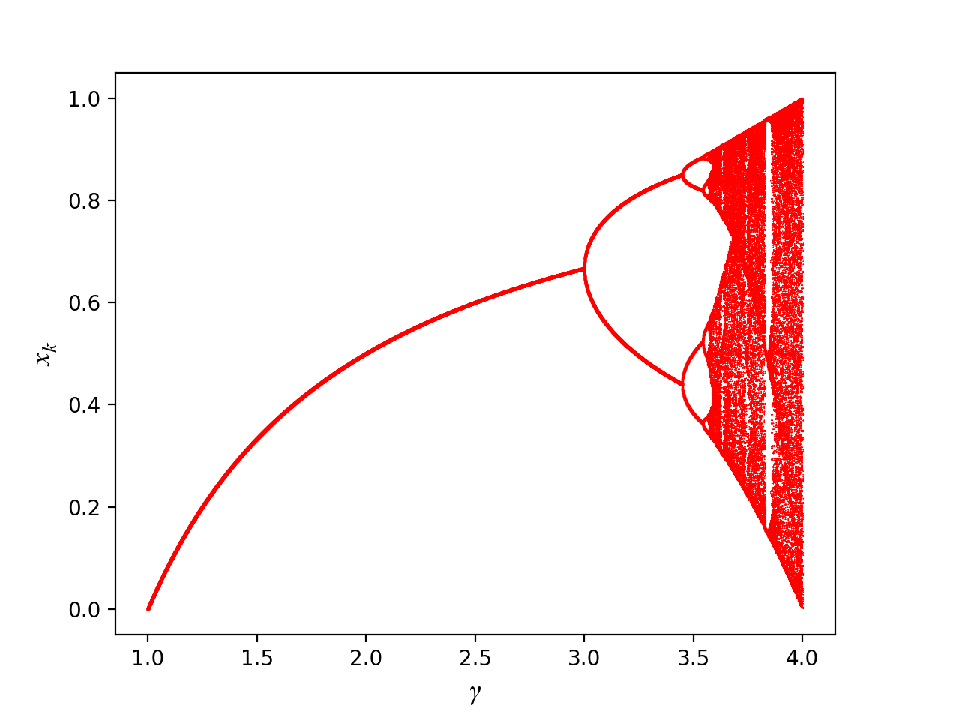}
        }
        \vspace{1em}
        \subfloat[$\delta=0.05$\label{fig:bif_compare:b}]{
            \includegraphics[width=0.696\columnwidth]{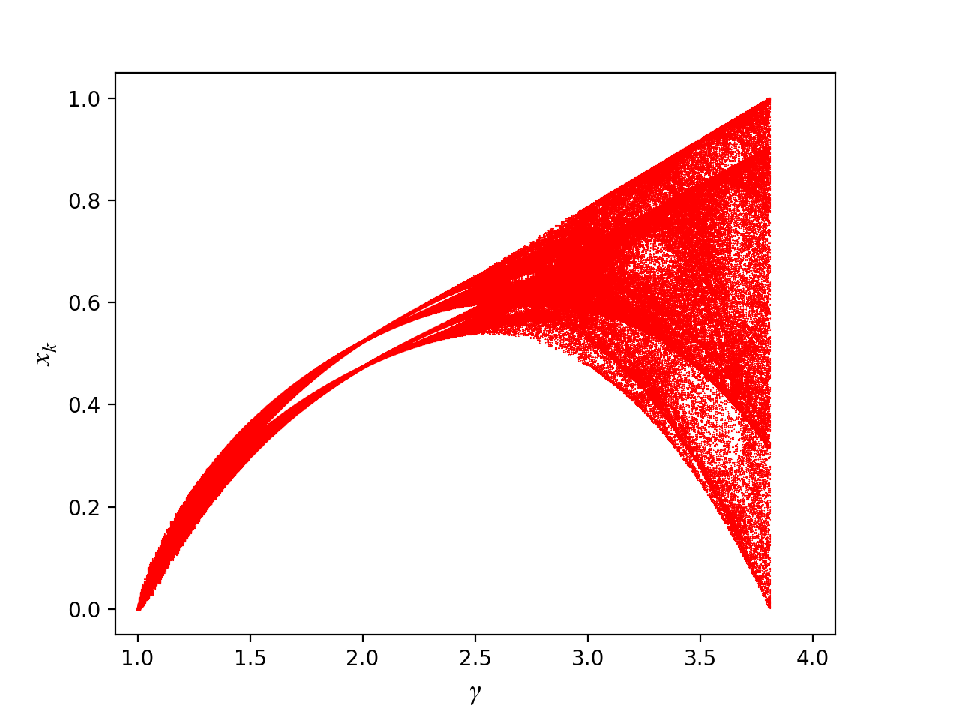}
        }
        \vspace{1em}
        \subfloat[$\delta=1/3$\label{fig:bif_compare:c}]{
            \includegraphics[width=0.696\columnwidth]{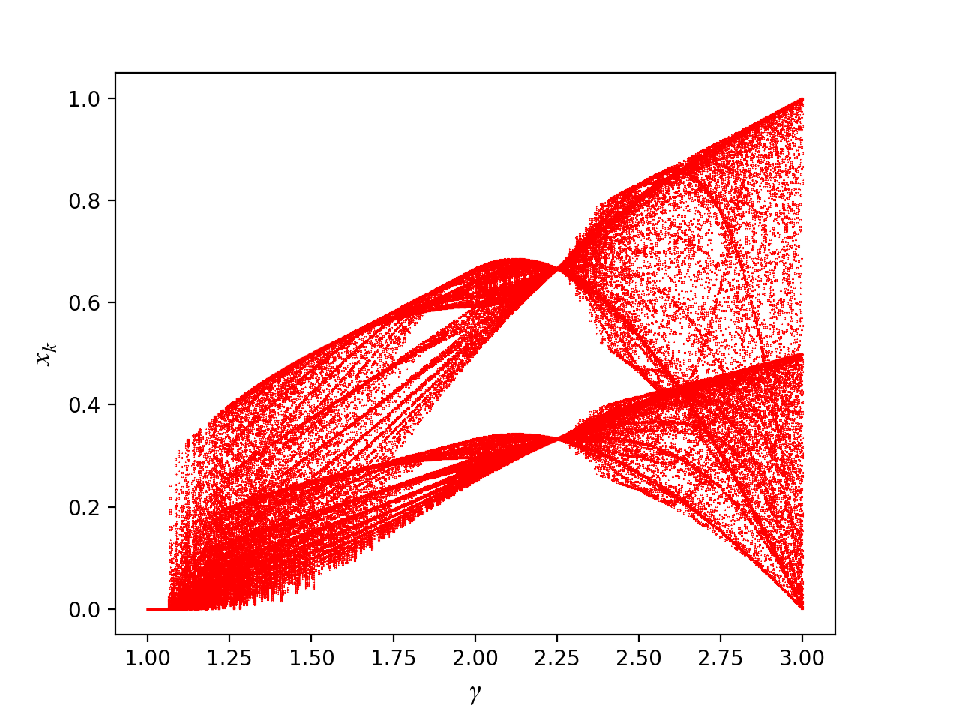}
        }
        \vspace{1em}
        \caption{\textbf{Bifurcation diagram of logistic IFS.} 
        Panel (a) shows the deterministic logistic map, corresponding to parameter gap $\delta=0$, while panels (b) and (c) correspond to $\delta=0.05$ and $\delta=1/3$, respectively.
        As shown in panel (b), even with a small $\delta$, the IFS typically exhibits noisy dynamics with blur in the bifurcation diagram in contrast to the deterministic case. In panel (c), at $\gamma=2.25$, the attractor consists of two points, demonstrating the emergence of a two-point toss-and-catch structure.
        }
        \label{fig:bif_compare}
    \end{figure}
\end{nolinenumbers}

\textbf{Figure~\ref{fig:bif_compare}} illustrates the bifurcation diagram with parameters $\alpha=\gamma (1-\delta)$ and $\beta=\gamma (1+\delta)$, where the parameter gap $\delta > 0$.
\textbf{Fig.~\ref{fig:bif_compare:a}} reproduces the bifurcation diagram of deterministic logistic map, whereas \textbf{Fig.~\ref{fig:bif_compare:b} and \ref{fig:bif_compare:c}} show the corresponding IFS dynamics with $\delta=0.05$ and $1/3$, respectively.
As shown in \textbf{Fig.~\ref{fig:bif_compare:b}}, even with small $\delta$, the IFS exhibits noisy dynamics with blur as opposed to the deterministic case.
In \textbf{Fig.~\ref{fig:bif_compare:c}}, we see the attractor consisting of two points at $\gamma=2.25$.

According to Abbasi \etal~(2018)~\cite{abbasi2018iterated}, even for relatively large values of \( \delta \), there exist pairs $(\alpha, \beta)$ for which the attractor of the IFS consists of exactly two points. These are examples of the toss-and-catch structure. Among the points in $\Lambda$, we call a point $B \in \Lambda \setminus \left(\Lambda_{g} \cup \Lambda_{h}\right)$ a {\bf bridging point}, that is, $B$ is not a periodic point of either $g$ or $h$, as shown in \eqref{eq:bridging}. In this paper, we focus on identifying finite invariant sets.

The return map corresponding to the 2-point toss-and-catch is shown in \textbf{Fig.~\ref{fig: 2points}}, while a schematic illustration of the relationship between the two points is provided in \textbf{Fig.~\ref{fig:2p-tac}}. We denote the period-$k$ points of $g$ (or $h$) by $g^{(k)}$ (or $h^{(k)}$) respectively.
As illustrated in {\bf Fig.~\ref{fig:2p-tac}}, the image of the fixed point $g^{(1)}$ by the function $h$ coincides with the fixed point $h^{(1)}$, and vice versa. That is, the state bounces between the fixed points of $g$ and $h$, like tossing a ball back and forth.

For the logistic map $f_\gamma(x)$, the nontrivial fixed point is given by $\dfrac{\gamma-1}{\gamma}$, and its preimages are $\left\{\dfrac{1}{\gamma}, \dfrac{\gamma-1}{\gamma}\right\}$. Using this property, the necessary and sufficient condition for the occurrence of a 2-point toss-and-catch can be expressed as
\[
    \alpha = \frac{\beta}{\beta-1}\quad (\beta \ne 1). \tag{Cond.~2}
    \label{eq: C2}
\]
As reported by Abbasi \etal~(2018)~\cite{abbasi2018iterated}, the expected Lyapunov exponent indicates that parameter variations can destabilize the toss-and-catch structure. Further details are provided in Sec.~\ref{sec:stability}. 

\begin{nolinenumbers}
    \begin{figure}[tb]
        \centering
        \subfloat[2-point toss-and-catch ($\alpha = 1.5~\text{and}~\ \beta=3.0$)\label{fig: 2points}]{
            \includegraphics[width=0.473\columnwidth]{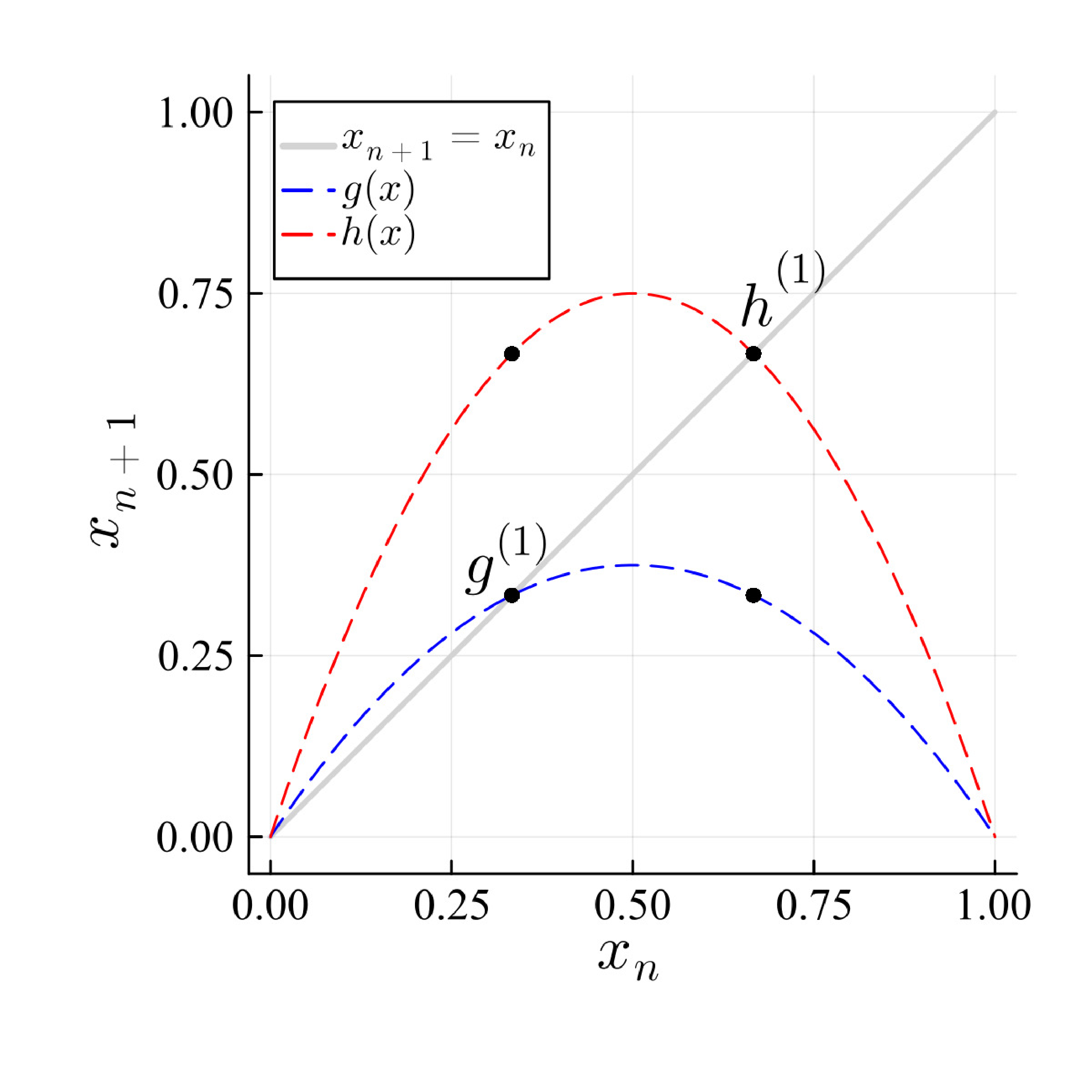}
        }\hfill
        \subfloat[3-point toss-and-catch ($\alpha \approx 2.3247~\text{and}~\ \beta\approx 3.0796$)\label{fig:3points}]{
            \includegraphics[width=0.473\columnwidth]{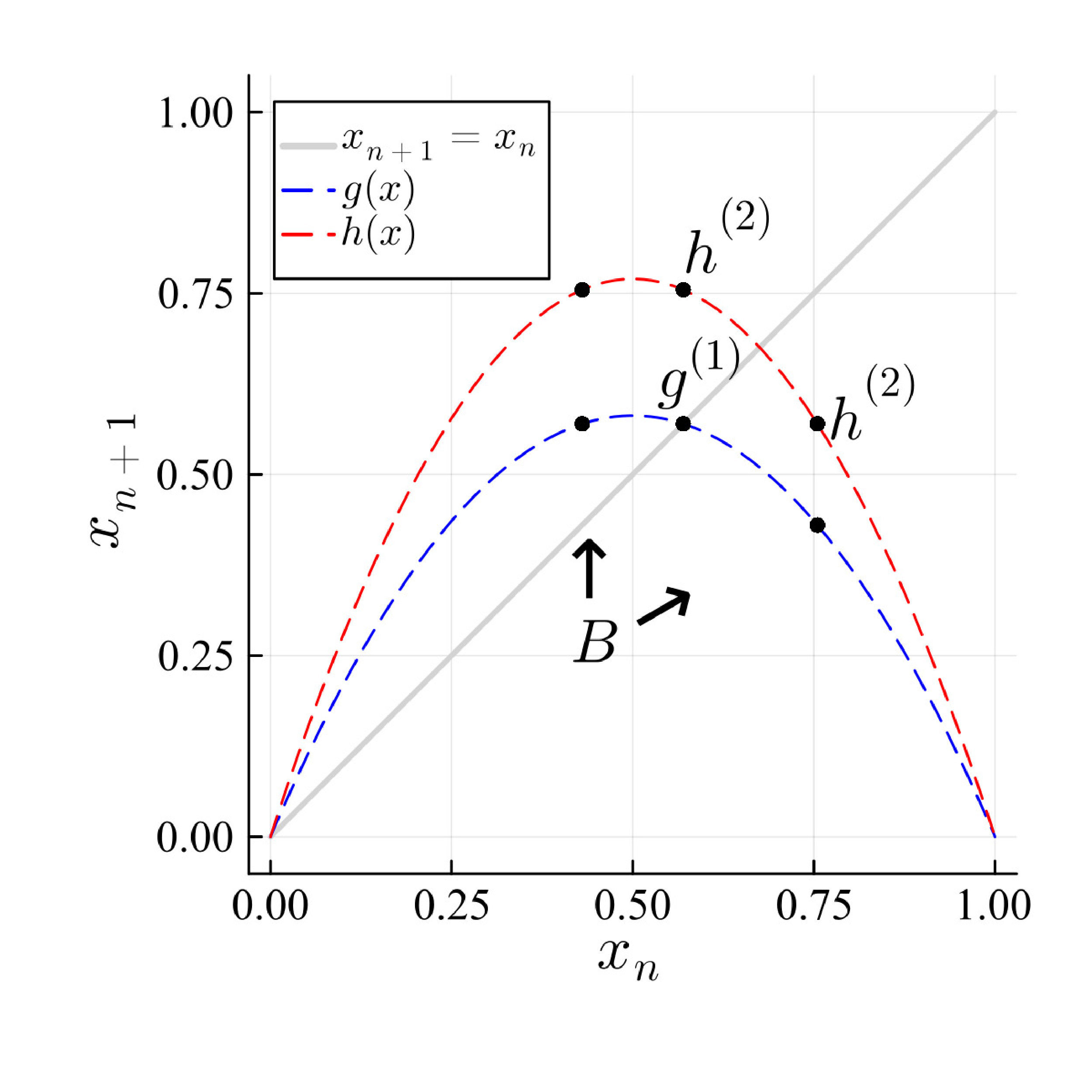}
        }\vspace{1em}
        \subfloat[5-point toss-and-catch ($\alpha \approx 2.1915~\text{and}~\ \beta\approx 3.3830$)\label{fig:5points}]{
            \includegraphics[width=0.473\columnwidth]{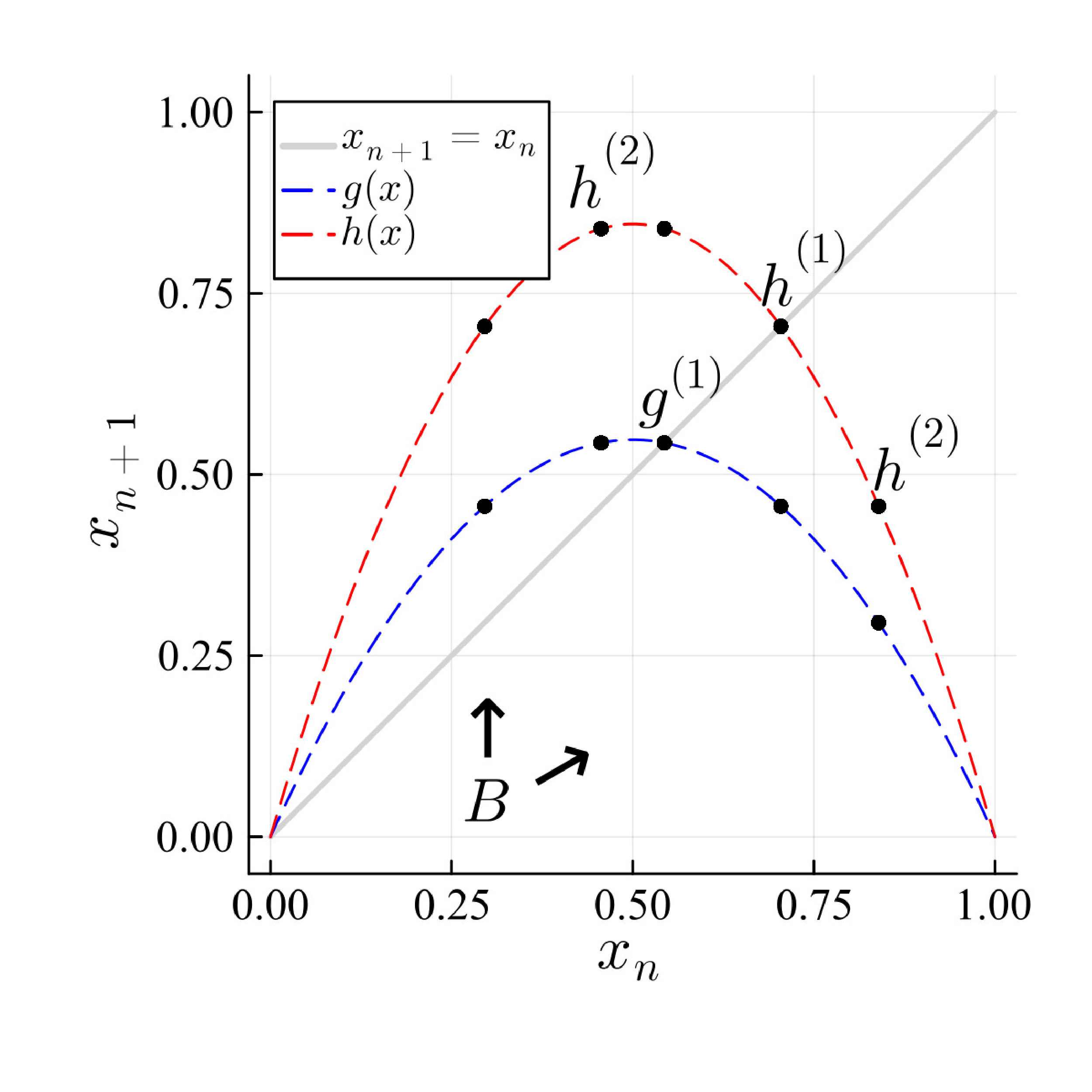}
        }
        \caption{\textbf{Return maps of 2-, 3-, and 5-point toss-and-catch of logistic IFS.} The blue and red dashed lines denote the graph of $g = f_\alpha$ and $h=f_\beta$, respectively. $g^{(k)}~(\text{or~} h^{(k)})$ denotes a period-$k$ point of the map $g~(\text{or~}h)$. $B$ denotes a bridging point that is not a periodic point of either $g$ or $h$.
        }
        \label{fig:network-structure}
    \end{figure}

    \begin{figure}[tb]
        \centering
        \subfloat[2-point toss-and-catch.
            \label{fig:2p-tac}
        ]{
            \includegraphics[width=0.47\linewidth]{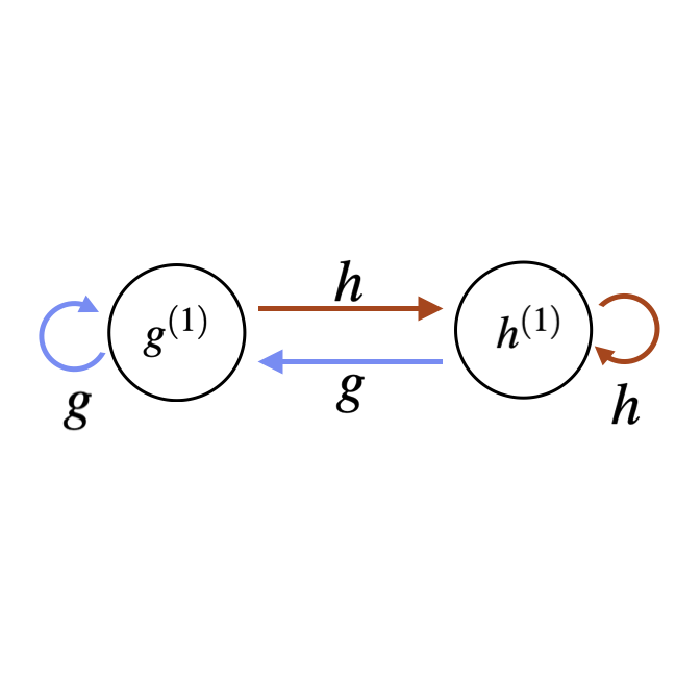}
        }\hfill
        \subfloat[3-point toss-and-catch.
            \label{fig:3p-tac}
        ]{
            \includegraphics[width=0.47\linewidth]{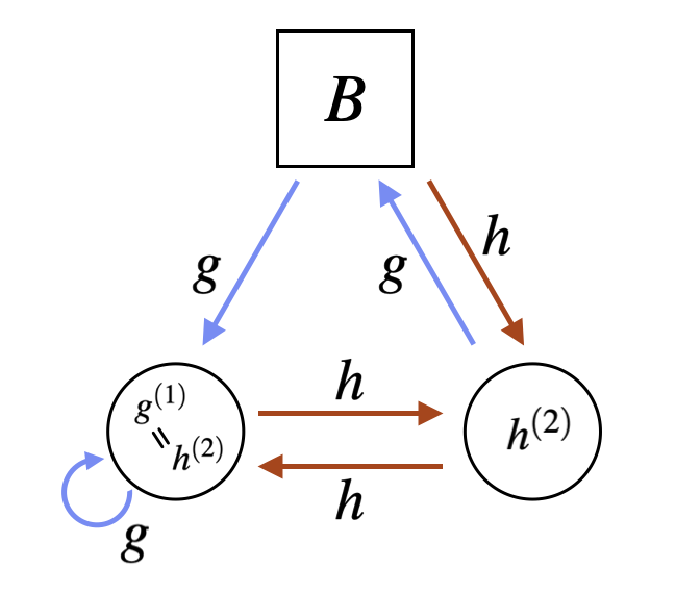}
        }
        \vspace{1em}
        \subfloat[5-point toss-and-catch.
            \label{fig:5p-tac}
        ]{
            \includegraphics[width=0.48\linewidth]{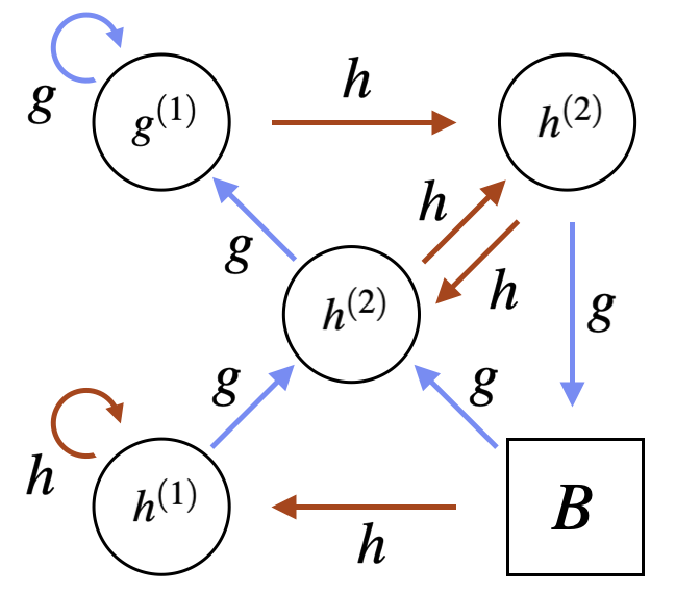}
        }

        \caption{\textbf{Structures of toss-and-catch dynamics in the logistic IFS.}
        Each panel illustrates the transitions between individual points in the toss-and-catch process. $g^{(k)}~(\text{or~} h^{(k)})$ denotes a period-$k$ point of map $g~(\text{or~}h)$. Panels (b) and (c) include bridging points $B$ that are not periodic under either $g$ or $h$. Panel (b) represents the simplest case of Eq.~\eqref{eq:bridging}.}
        \label{fig:tac-structure}
    \end{figure}
\end{nolinenumbers}

For a different pair of $(\alpha, \beta)$, the system exhibits a 3-point toss-and-catch structure.
\textbf{Figures~\ref{fig:3points},~\ref{fig:3p-tac}} present the corresponding return map and graphical structure of this 3-point toss-and-catch. In this case, the attractor consists of a stable fixed point of $g$ and points on a stable period-2 orbit of $h$.
Notably, the components of the toss-and-catch are not limited to invariant sets of the logistic map but also include the bridging point $B$. It bridges the periodic points of $g$ and $h$, and complements the network.
The period-2 orbit of $h$ is given by
\[
    \left\{ \frac{\beta+1 - \sqrt{(\beta+1)(\beta-3)}}{2\beta},\ \frac{\beta+1 + \sqrt{(\beta+1)(\beta-3)}}{2\beta} \right\}.
\]
The following condition (Cond.~3) is derived from the geometric relationships among the points in Fig.~\ref{fig:3p-tac}:
\[
    \begin{cases}
        \begin{array}{rcccl}
            & \frac{\beta+1 - \sqrt{(\beta+1)(\beta-3)}}{2\beta} &=& \frac{\alpha-1}{\alpha}, & \quad \text{\small(Cond.~3a)} \\
            & h\left(\frac{1}{\alpha}\right) &=& \frac{\beta+1 + \sqrt{(\beta+1)(\beta-3)}}{2\beta}, & \quad \text{\small(Cond.~3b)} \\
            & g\left(\frac{\beta+1 + \sqrt{(\beta+1)(\beta-3)}}{2\beta}\right) &=& \frac{1}{\alpha}. & \quad \text{\small(Cond.~3c)}
        \end{array}
    \end{cases}
\]
By solving the conditions (Cond.~3a), (Cond.~3b), and (Cond.~3c), we obtain the parameter values for which 3-point toss-and-catch occurs:
\begin{align*}
    \begin{cases}
        \alpha &= \dfrac{1}{3} \left(3 + \sqrt[3]{\tfrac{1}{2}\left(27 - 3\sqrt{69}\right)}
            + \sqrt[3]{\tfrac{1}{2}\left(27 + 3\sqrt{69}\right)}
        \right), \\
        \beta &= \dfrac{1}{3} \left(
            2 + \sqrt[3]{\tfrac{1}{2}\left(97 - 3\sqrt{69}\right)}
            + \sqrt[3]{\tfrac{1}{2}\left(97 + 3\sqrt{69}\right)}
        \right).
    \end{cases}
\end{align*}
Numerically, these values correspond to $(\alpha, \beta) \approx (2.3247 , 3.0796)$.

\textbf{Figures~\ref{fig:5points}} and \textbf{\ref{fig:5p-tac}} show the return map and the graph of 5-point toss-and-catch structure.
In the 5-point case, the attractor consists of a stable fixed point of $g$, an unstable fixed point and a stable period-2 orbit of $h$, and a bridging point B.

The following condition (Cond.~5) arises from the relationships among the points in Fig.~\ref{fig:5p-tac}:
\[
    \begin{cases}
        \begin{array}{rcccl}
            & h\left(\frac{\alpha - 1}{\alpha}\right) &=& \frac{\beta+1 + \sqrt{(\beta+1)(\beta-3)}}{2\beta}, & \quad \text{\small(Cond.~5a)} \\
            & g\left(\frac{\beta+1 - \sqrt{(\beta+1)(\beta-3)}}{2\beta}\right) &=& \frac{\alpha-1}{\alpha}, & \quad \text{\small(Cond.~5b)} \\
            & g\left(\frac{\beta-1}{\beta}\right) &=& \frac{\beta+1 - \sqrt{(\beta+1)(\beta-3)}}{2\beta}, & \quad \text{\small(Cond.~5c)} \\
            & g\left(\frac{1}{\beta}\right) &=& \frac{\beta+1 - \sqrt{(\beta+1)(\beta-3)}}{2\beta}, & \quad \text{\small(Cond.~5d)} \\
            &    g\left(\frac{\beta+1 + \sqrt{(\beta+1)(\beta-3)}}{2\beta}\right) &=& \frac{1}{\beta}. & \quad \text{\small(Cond.~5e)}
        \end{array}
    \end{cases}
\]
By solving the conditions (Cond.~5a)--(Cond.~5e), we obtain the following parameter values at which 5-point toss-and-catch occurs:
\begin{align*}
    \begin{cases}
        \alpha &= 1 + \dfrac{1}{6}\left(\sqrt[3]{54 - 6\sqrt{33}} + \sqrt[3]{54 + 6\sqrt{33}}\right),\\
        \\
        \beta &= 1 + \dfrac{1}{3}\left(\sqrt[3]{54 - 6\sqrt{33}} + \sqrt[3]{54 + 6\sqrt{33}}\right).
    \end{cases}
\end{align*}
These correspond approximately to $(\alpha, \beta) \approx (2.1915, 3.3830)$.

To systematically locate parameter sets exhibiting toss-and-catch dynamics, we quantify the “size” of the IFS attractor.
In \textbf{Fig.~\ref{fig: phase-diagram}}, this is measured by the minimum number of \(\varepsilon\)-balls required to cover the attractor (see the supplementary material for details).  Note that this method does not capture unstable toss-and-catch structures. In \textbf{Fig.~\ref{fig: phase-diagram:a}}, blue regions indicate parameter sets whose attractors can be covered by only a few \(\varepsilon\)-balls (here, \(\varepsilon=10^{-6}\)). The diagonal \(\beta=\alpha\), corresponding to the deterministic logistic map, remains blue up to the onset of chaos. Parameter values with superstable periodic orbits (e.g., \(\gamma=2, 1+\sqrt{5}\)) also appear blue due to rapid convergence.

A distinct curve emerges, marking regions where the IFS attractor is sparse relative to the surrounding dense regions. The parameter curves determined by conditions (Cond.~2), (Cond.~3), and (Cond.~5) are shown in {\bf Fig.~\ref{fig: phase-diagram:b}}.
\allred 
We observe a toss-and-catch-like structure near the curve defined by \eqref{eq: C2}. Although not a precise toss-and-catch, this structure is relatively coarse-grained: in parameter regions near the curve defined by (Cond.~2), the attractor remains sufficiently sparse to be covered by two $\varepsilon$-balls for some small $\varepsilon$.
This coarse-grained persistence is observed near exact toss-and-catch parameter sets for which the corresponding exact finite invariant set has a negative expected Lyapunov exponent, as discussed in Sec.~~\ref{sec:stability}.
This observation naturally raises the question of why such structures are sometimes not clearly observed in numerical simulations.

If a parameter set $(\alpha,\beta)$ satisfies conditions (Cond.~2), (Cond.~3), or (Cond.~5), the corresponding finite invariant set exists regardless of the value of the switching probability $p$.
In contrast, the invariant (stationary) distribution on $\Lambda$ generally depends on $p$, because the transition probabilities between points of $\Lambda$ are themselves $p$-dependent.
Consequently, while the geometric existence of toss-and-catch structures is independent of $p$, their observability and linear stability depend on the stationary measure induced by the random switching process.

To clarify this distinction between geometric existence and probabilistic stability, we next examine the stationary measures and expected Lyapunov exponents associated with toss-and-catch dynamics.

\begin{nolinenumbers}
    \begin{figure}[]
        \centering
        \subfloat[Attractor ``size'' for parameter $(\alpha, \beta)$\label{fig: phase-diagram:a}]{
            \includegraphics[width=0.87\columnwidth]{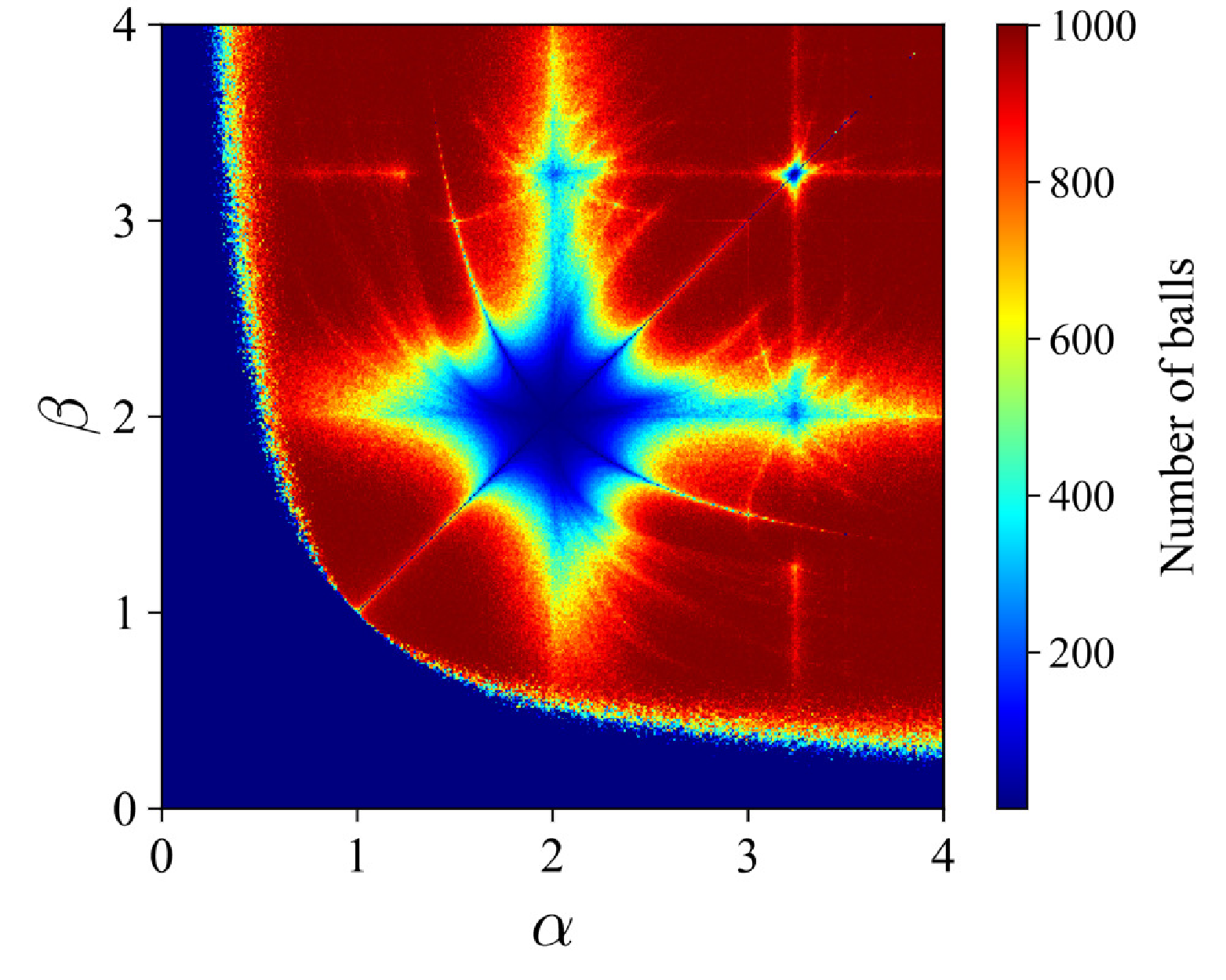}
        }\hfill
        \subfloat[Parameter relations of (Cond.~2), (Cond.~3), and (Cond.~5)\label{fig: phase-diagram:b}]{\includegraphics[width=0.87\columnwidth]{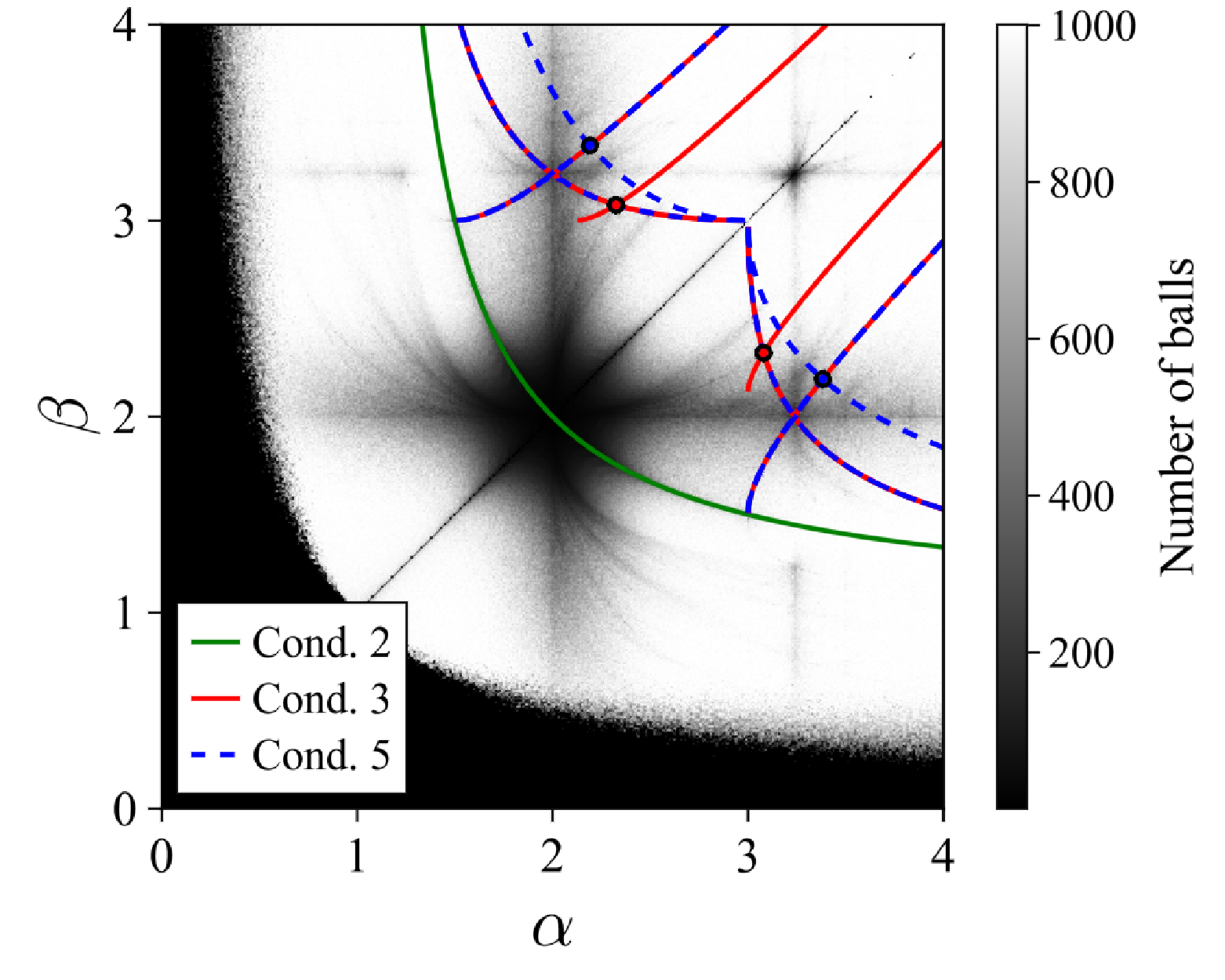}}
        \caption{\textbf{Approximation of attractor ``size'' for each parameter set of the logistic IFS.}
            For each parameter set $(\alpha, \beta)$, the logistic IFS is iterated for 1000 steps after discarding transients.
            {\bf Panel~(a)} shows the minimum number of $\varepsilon$-balls ($\varepsilon = 10^{-6}$) required to cover the attractor, the result represented using a color scale.
            The region highlighted in blue corresponds to parameters for which the attractor can be covered with fewer $\varepsilon$-balls.
            {\bf Panel~(b)} illustrates the sets of parameters satisfying each condition of (Cond.~2), (Cond.~3), and (Cond.~5). A 2-point toss-and-catch occurs along the curve defined by condition (Cond.~2), while 3- and 5-point toss-and-catch appear at the intersections of the curves defined by conditions (Cond.~3) and (Cond.~5), respectively. These two points are highlighted with dots.
        }
        \label{fig: phase-diagram}
    \end{figure}
\end{nolinenumbers}

\section{Stability and Stationary Measure}\label{sec:stability}
In Secs. II and IV, we fixed the switching probability at $p = 0.5$ to focus on the geometric structure of toss-and-catch dynamics. 
In contrast, in this section, we regard $p$ as the primary control parameter and examine how the stationary measure and stability properties depend on the switching process.

To understand why toss-and-catch structures are sometimes difficult to observe numerically, we examine the stationary measures and expected Lyapunov exponents associated with the induced Markov dynamics on finite invariant sets.
Restricted to a finite invariant set $\Lambda$, the IFS defines a finite-state Markov chain whose transition probabilities are determined by the switching probability $p$. 
Let $P(p)$ denote the row-stochastic transition matrix on $\Lambda$, where $P_{ij}(p)$ represents the transition probability from the $i$-th point to the $j$-th point. 
We write stationary distributions as row vectors. 
The stationary distribution $\pi$ satisfies
\[
\pi = \pi P(p).
\]
Although the existence of toss-and-catch structures is independent of $p$, their stability depends on the stationary measure induced by $P(p)$ and on the derivatives of the selected maps along $\Lambda$.
Numerical experiments indicate that for certain values of $p$, the invariant set becomes unstable and the nearby trajectory escapes to a larger invariant set.

For the 2-point toss-and-catch, the induced Markov chain has a transition matrix
\[
P=
\begin{pmatrix}
p & 1-p\\
p & 1-p
\end{pmatrix},
\]
with stationary distribution $\pi=(p,1-p)$.
The expected Lyapunov exponent $\lambda_p$ is
\[
E[\lambda_p]
=
p\ln|2-\alpha|
+
(1-p)\ln\left|\frac{\alpha-2}{\alpha-1}\right|.
\]
This demonstrates that while the invariant set itself exists independently of $p$, its stability depends continuously on the switching probability.

For the 3-point toss-and-catch, the stationary measure becomes asymmetric 
because transitions through the bridging point are probabilistically asymmetric.
This illustrates how switching-induced invariant structures influence not only geometry but also probabilistic stability.
These results highlight the distinction between geometric existence and dynamical observability of toss-and-catch structures.
Additional calculations for higher-point toss-and-catch structures and logistic-tent IFS cases are provided in the Supplementary Material.
\allblack

\section{Logistic-Tent IFS}\label{sec:logistic-tent}
\begin{figure}
    \centering
    \includegraphics[width=0.98\linewidth]{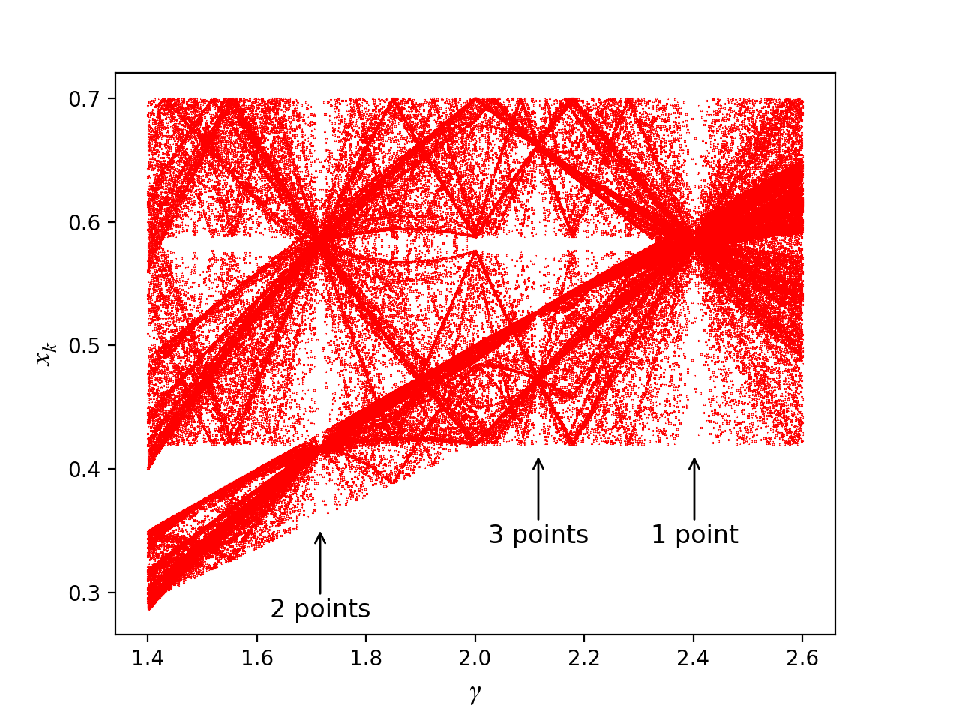}
    \caption{\textbf{Bifurcation diagram of logistic-tent IFS.} For $\mu = 1.4$ and $\gamma \in [1.4, 2.6]$, we obtained the bifurcation diagram of the logistic-tent IFS. The $n$-point toss-and-catch occurs at $\gamma = 2.4$~($n=1$), $\gamma \approx 1.7143$~($n=2$), and $\gamma \approx 2.1143$~($n=3$). The return maps for these cases are shown in \textbf{Fig.~\ref{fig:returnmap-lt}}. }
    \label{fig:bif-lgs-tent}
\end{figure}

    \begin{figure}
        \centering
        \subfloat[1-point toss-and-catch ($\mu= 1.4,~\gamma = 2.4$). \label{fig: return1_lt}]{
            \includegraphics[width=0.478\columnwidth]{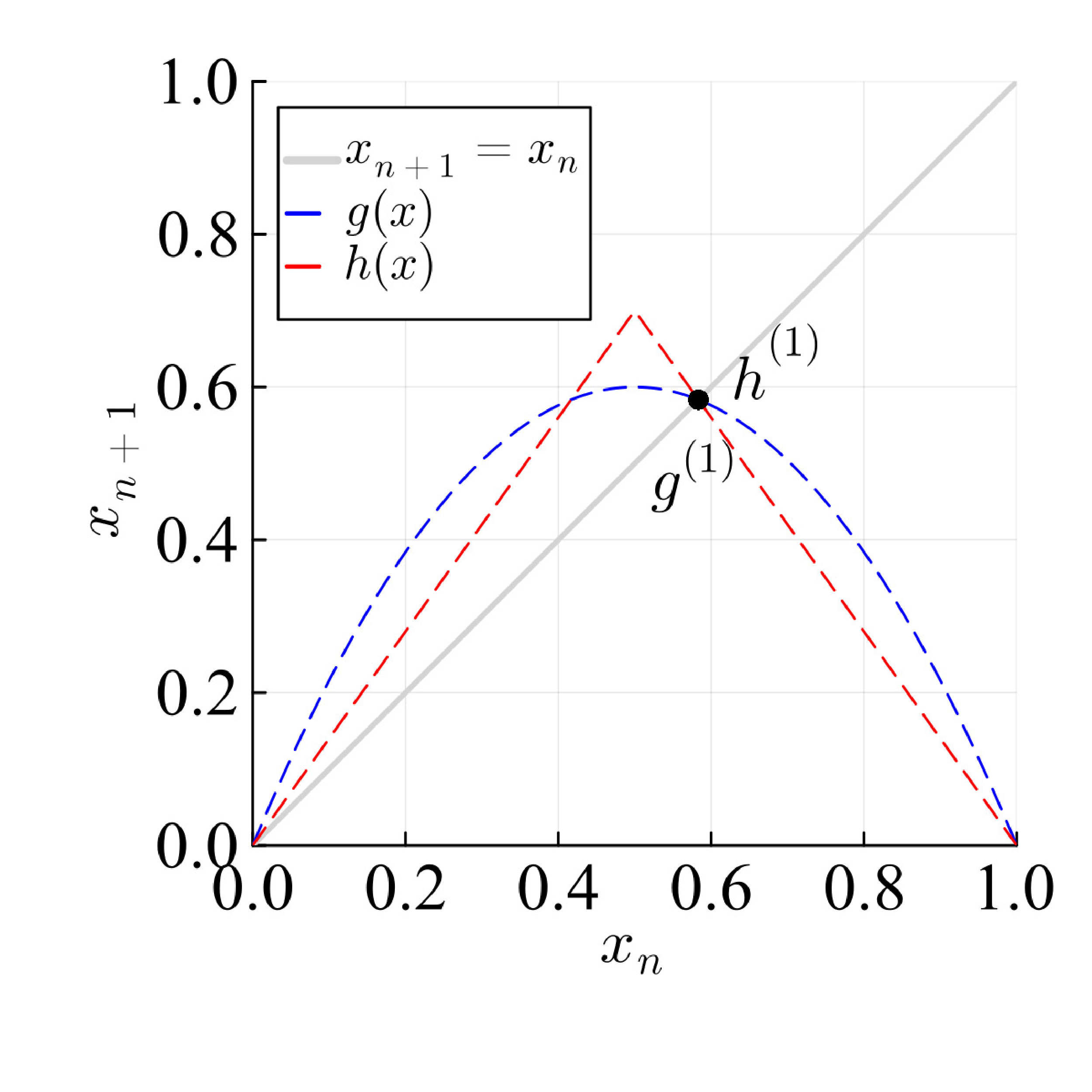}
        }\hfill
        \subfloat[2-point toss-and-catch ($\mu = 1.4,~\gamma \approx 1.7143$).\label{fig: return2_lt}]{
            \includegraphics[width=0.478\columnwidth]{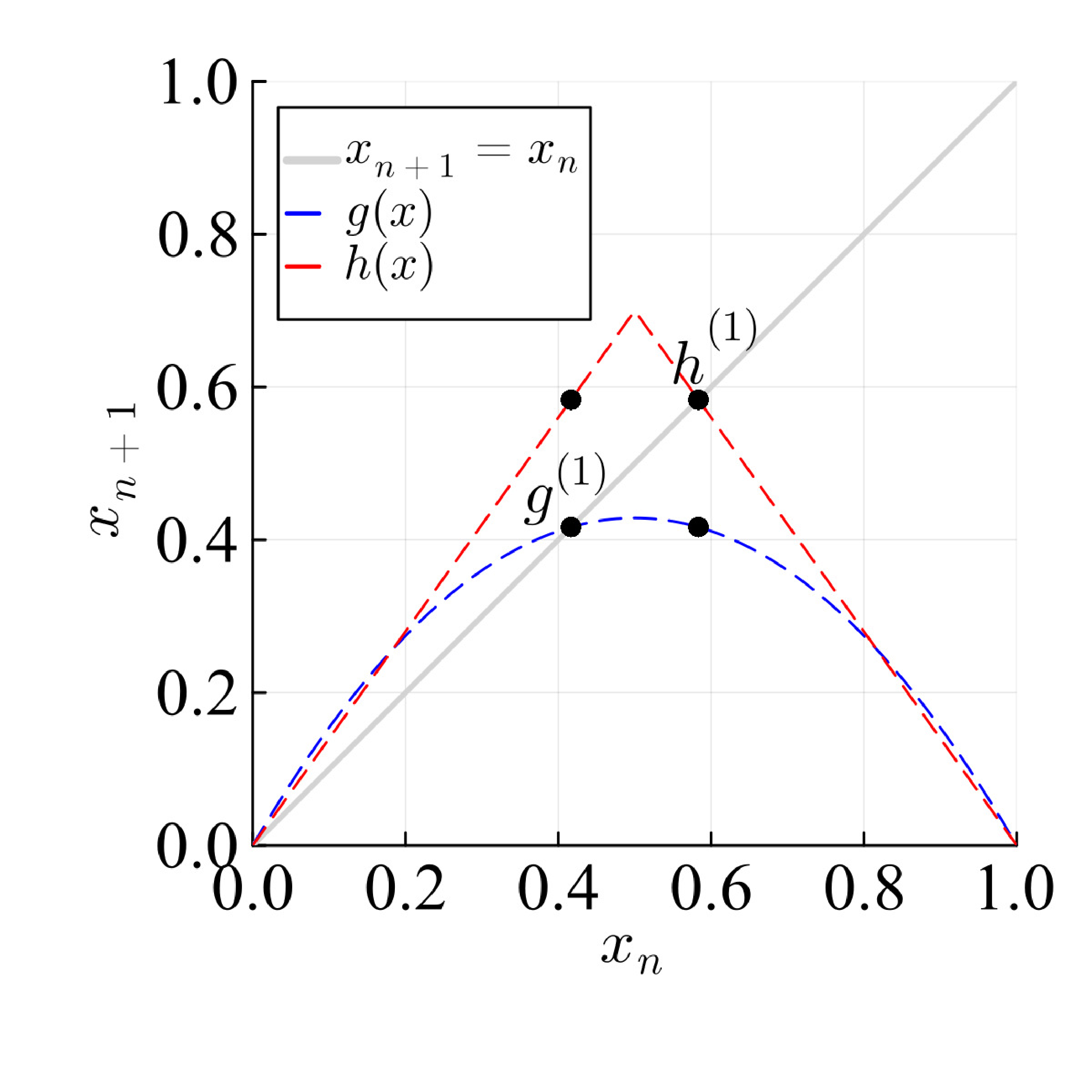}
        }\vspace{1em}
        \subfloat[3-point toss-and-catch ($\mu= 1.4,~\gamma \approx 2.1143$).\label{fig: return3_lt}]{
            \includegraphics[width=0.478\columnwidth]{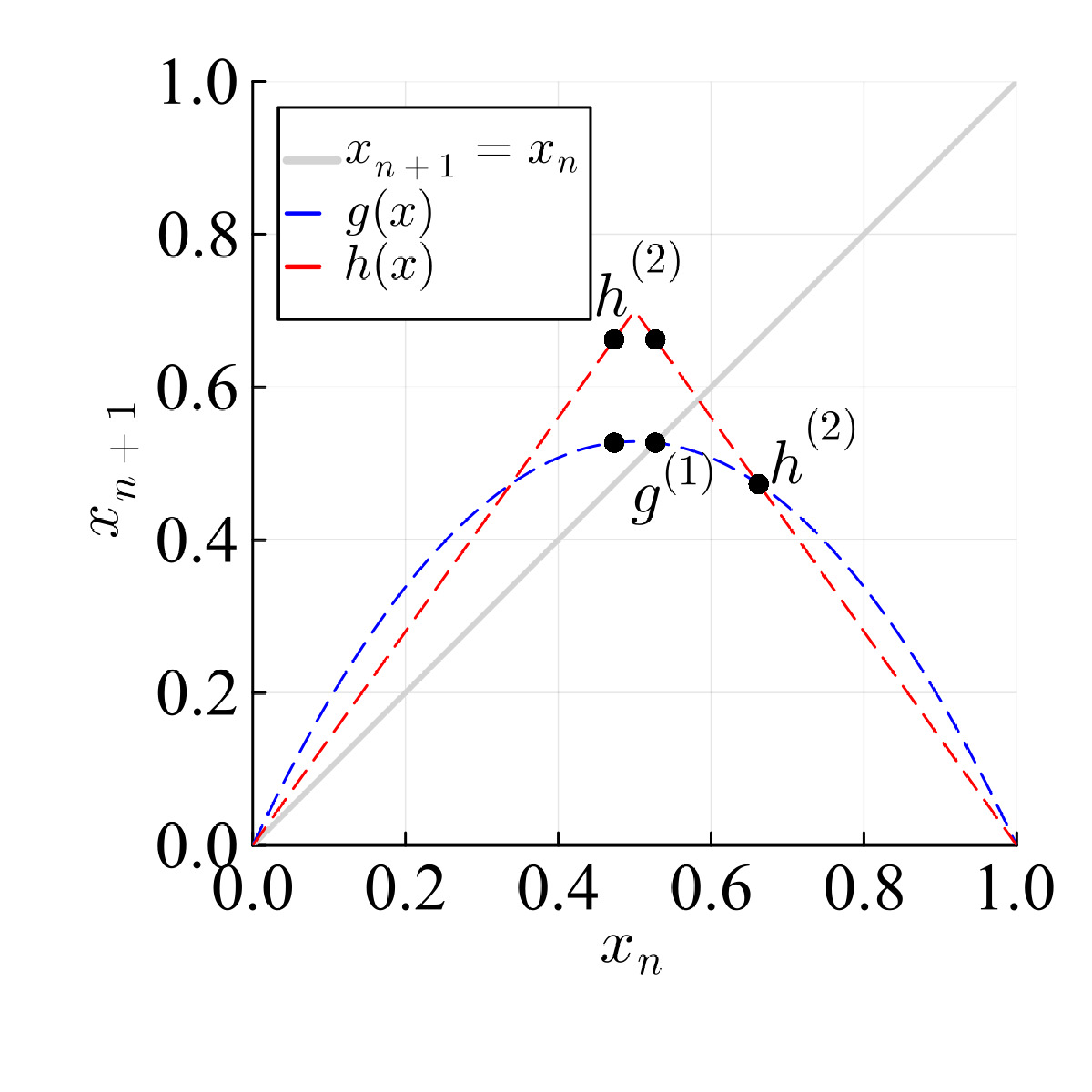}
        }
        \caption{\textbf{Return maps of 1-, 2-, and 3-point toss-and-catch of logistic-tent IFS.} For the 2-point toss-and-catch case, the structure resembles that of the logistic IFS. However, the 3-point toss-and-catch structure differs from the logistic IFS because the node $B$ is absent.
        }\label{fig:returnmap-lt}
    \end{figure}

In this section, we consider a random dynamical system that alternates between the logistic and tent maps. For $\mu \in [0,2]$, $\gamma \in [0,4]$, the logistic-tent IFS is defined as:
\begin{align*}
    F(x)=
    \begin{cases}
        \kern-0.3em
        \begin{array}{lll}
            g(x) &= \gamma x(1-x), & \kern-1.5em \text{with }p,\\
            h(x) &=
            \begin{cases}
                \begin{array}{rl}
                    \mu x, & x \in [0,~0.5)\\
                    \kern-0.8em -\mu (x-1), & x \in [0.5,~1]
                \end{array}
            \end{cases} & \kern-1.5em \text{with }1-p,
        \end{array}
    \end{cases}
\end{align*}
where we set $p=0.5$ in the numerical examples of this section, and the maps are selected i.i.d. at each step.
We investigate whether the toss-and-catch dynamics observed in the logistic IFS also occur in this system. If $\mu = 1$, the left half of the tent map becomes the identity, producing a trivial toss-and-catch. 
For the nontrivial toss-and-catch structures discussed below, we focus on the case $\mu>1$, where the tent map has a nonzero fixed point on its right branch. The logistic-tent IFS has nontrivial intersection points, which do not appear in the logistic IFS. These intersections affect the structure of the 3-point toss-and-catch.

\textbf{Figure~\ref{fig:bif-lgs-tent}} shows the bifurcation diagram of the logistic-tent IFS obtained by fixing $\mu = 1.4$ and varying $\gamma$ over the interval $[1.4, 2.6]$. In this range, toss-and-catch structures with one, two, and three points are observed.
\textbf{Figure~\ref{fig:returnmap-lt}} illustrates the corresponding return maps.

A one-point toss-and-catch occurs when the parameters of the two maps satisfy the condition:
\[
    \gamma = 1+\mu \tag{Cond.~1$'$}.
\]
These points correspond to the coincidence of the fixed points of the logistic and tent maps (\textbf{Fig.~\ref{fig: return1_lt}}). In the logistic IFS, this occurs when $g$ and $h$ are identical.

A 2-point toss-and-catch occurs when the parameters satisfy the condition:
\[
    \gamma = \frac{1+\mu}{\mu} \tag{Cond.~2$'$}.
\]
The two points are the nontrivial fixed point of $g$ and the nontrivial fixed point of $h$, respectively:
\[
\left\{1-\frac{1}{\gamma}, \frac{\mu}{1+\mu}\right\}.
\]
This structure is the same as the 2-point toss-and-catch of logistic IFS (\textbf{Fig.~\ref{fig:2p-tac}}).

\begin{figure}
    \centering
    \begin{minipage}{0.478\columnwidth}
        \includegraphics[width=\linewidth]{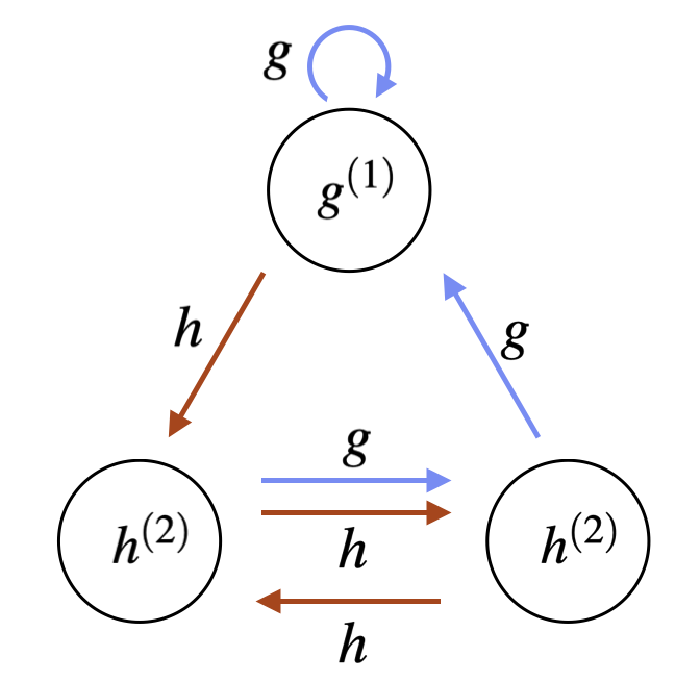}
    \end{minipage}
    \caption{\textbf{Structure of 3-point toss-and-catch in logistic-tent IFS.} $g^{(k)}~(\text{or~} h^{(k)})$ denotes a period-$k$ point of the map $g~(\text{or~}h)$. Unlike the 3-point toss-and-catch structure in the logistic IFS, no bridging point exists here. Further discussion is provided in the supplementary material.}
    \label{fig:network-lgs-tent}
\end{figure}

    \begin{figure}
        \centering
        \includegraphics[width=\linewidth]{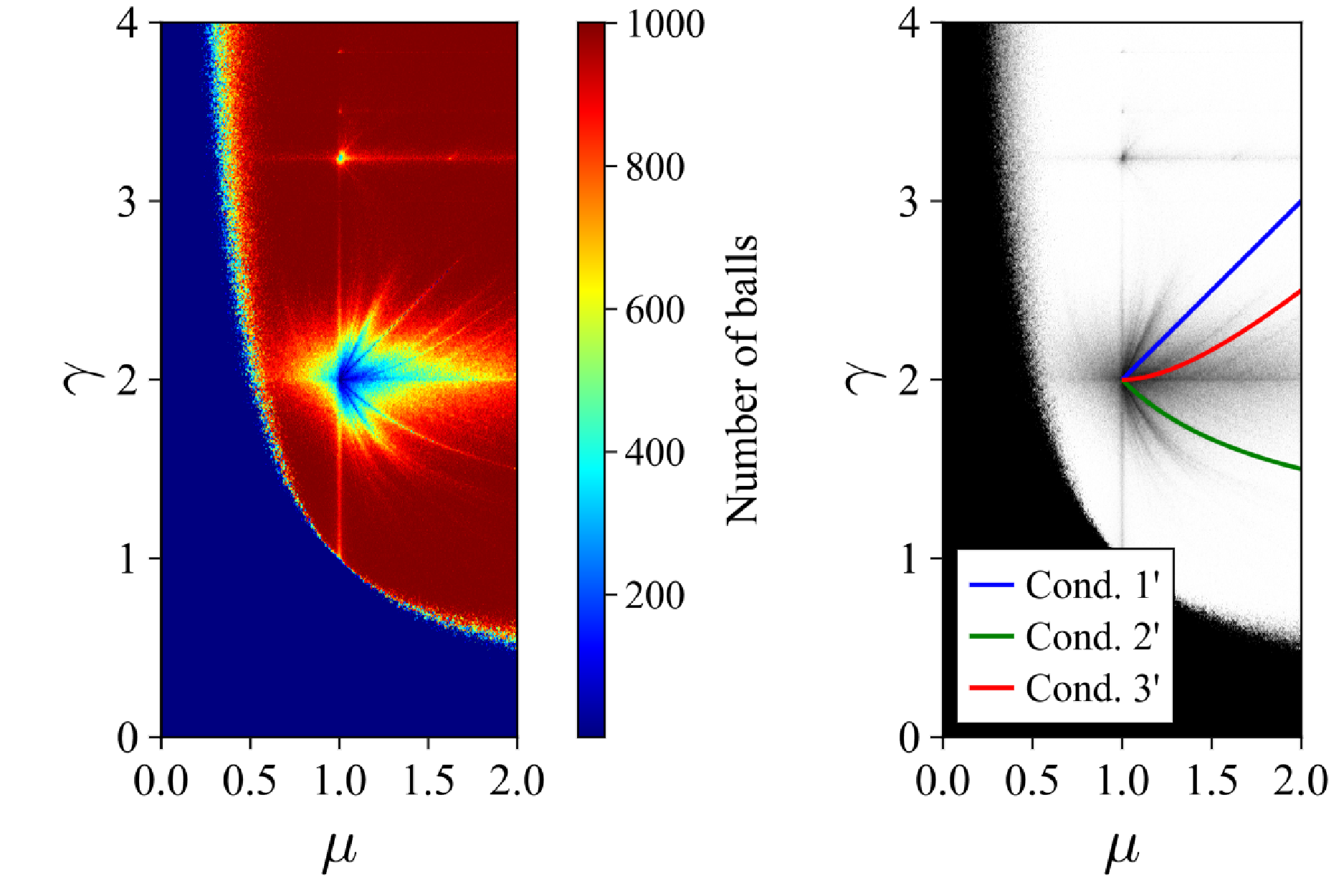}
        \caption{\textbf{Approximation of attractor ``size'' for each parameter of the logistic-tent IFS.}
            For each parameter pair $(\mu, \gamma)$, the logistic-tent IFS is iterated for 1000 steps after discarding transients.
            The left panel shows the number of $\varepsilon$-balls ($\varepsilon = 10^{-6}$) required to cover the attractor, displayed using a color scale.
            The region highlighted in blue corresponds to parameters for which the attractor can be covered with fewer $\varepsilon$-balls.
            The right panel illustrates parameter sets satisfying conditions (Cond.~1$'$), (Cond.~2$'$), and (Cond.~3$'$), corresponding to 1-, 2-, and 3-point toss-and-catch, respectively.
        }
        \label{fig:phase-diagram-lt}
    \end{figure}

A 3-point toss-and-catch structure occurs when the parameters satisfy the condition:
\[
    \gamma = \frac{1 + \mu^2}{\mu} \tag{Cond.~3$'$}.
\]
This set consists of the nontrivial fixed point of $g$ and the two points on a period-2 orbit of $h$:
\[
\left\{
\frac{\mu}{1+\mu^2},\,
1-\frac{1}{\gamma},\,
\frac{\mu^2}{1+\mu^2}
\right\}.
\]
Unlike the logistic IFS, this system consists only of periodic points of $g$ and $h$, as illustrated in {\bf Fig.~\ref{fig:network-lgs-tent}}. This behavior arises because $g$ and $h$ possess an intersection other than $\{0,\ 1\}$. In the logistic IFS, a contradiction appears in the 3-point case unless a bridging point associated with a preimage is introduced. In the logistic-tent IFS, this contradiction is avoided because the two maps share a nontrivial intersection point.
Further details are provided in the supplementary material.

{\bf Figure~\ref{fig:phase-diagram-lt}} illustrates an approximation of attractor ``size'' using $\varepsilon$-balls ($\varepsilon=10^{-6}$).
Similar to the logistic IFS, a distinct curve emerges, marking regions where the IFS attractor is sparse relative to the surrounding dense regions. The curves defined by the conditions (Cond.~1$'$), (Cond.~2$'$), and (Cond.~3$'$) are shown in the right panel. If the parameter set $(\mu, \gamma)$ satisfies conditions (Cond.~1$'$), (Cond.~2$'$), or (Cond.~3$'$), the corresponding finite invariant set exists regardless of the value of the switching probability $p$.
We observe a toss-and-catch-like structure near the curves, similar to that of the logistic IFS.
This observation suggests that stable toss-and-catch structures possess a certain coarse-grained robustness in the logistic-tent IFS: although the exact finite invariant set exists only along specific parameter curves, nearby parameter regions may still exhibit sparse and visually recognizable toss-and-catch-like dynamics when the corresponding expected Lyapunov exponent remains negative, as discussed in Sec.~\ref{sec:stability}.

\section{Concluding Remarks}
\subsection{Summary}
We have investigated IFS that randomly alternate between two distinct maps, $g$ and $h$. We analyze two representative cases:
$g$ and $h$ are both logistic maps (logistic IFS), and $g$ is a logistic map and $h$ is a tent map (logistic-tent IFS). Our focus is on toss-and-catch structures, which are finite invariant sets of the IFS.
In both the logistic IFS and logistic-tent IFS, we observe finite invariant toss-and-catch structures. These include structures formed by transitions between fixed points and structures in which fixed points are connected to period-2 orbits.
Notably, we identify cases in which the invariant set of the IFS $\Lambda$ contains bridging points, that is, points belonging to neither $\Lambda_g$ nor $\Lambda_h$.
These findings extend prior work on IFS dynamics.
We also observe that toss-and-catch structures exist along specific parameter curves, with the attractor remaining relatively sparse in neighboring parameter regions.

\allred
\subsection{Discussion}
From the perspective of emergent invariant structures in switching dynamical systems, the logistic IFS and logistic-tent IFS may be regarded as analytically tractable low-dimensional realizations of two representative mechanisms for finite invariant sets in one-dimensional random dynamical systems.
In the logistic IFS, the 3- and 5-point toss-and-catch structures require bridging points because the relevant periodic structures of the constituent maps cannot be connected solely through direct intersections.
In contrast, the logistic-tent IFS exhibits finite invariant sets arising from nontrivial intersections among the maps themselves.

The present systems provide analytically tractable low-dimensional models in which the closure relations underlying finite invariant structures can be identified explicitly. 
In higher-dimensional systems, similar mechanisms are expected to exist, but invariant objects are typically manifolds or other complicated invariant sets rather than finite collections of points. 
Consequently, bridging points may generalize to higher-dimensional bridging structures that connect invariant sets of different maps, while the associated geometric and stability properties become substantially more complex. 
Such structures may include higher-dimensional connecting invariant sets or fractal heteroclinic-like sets linking different invariant structures. 
This viewpoint is reminiscent of heteroclinic structures observed in heterochaotic systems~\cite{saiki_piecewise_2021,saiki23}, where fractal connecting sets arise between invariant sets associated with different unstable dimensions.
In some heterochaotic baker-type systems, the expanding dynamics generating symbolic itineraries may also be interpreted as inducing deterministic switching between distinct local dynamical behaviors, leading to fractal connecting structures.
These perspectives suggest that toss-and-catch structures are not merely isolated examples but represent a more general mechanism in switching dynamical systems.

Many studies of switching or regime-switching dynamical systems focus on the effects of switching on stability, synchronization, intermittency, or statistical properties of trajectories. 
From this perspective, the present results highlight a different aspect of switching dynamics: the emergence of invariant structures that do not exist in the constituent systems individually.
In particular, the bridging points identified in the logistic IFS are not periodic points of either map alone, but become essential components of the invariant set only through random switching between maps. 
This suggests that switching can generate qualitatively new invariant geometries rather than merely perturbing the dynamics of existing invariant sets.
The present study provides a simple, low-dimensional example that illustrates how regime switching may create emergent invariant structures through interactions among distinct dynamical rules.
\allblack

\section{Supplementary Material}
The supplementary material provides detailed calculations of transition matrices, stationary distributions, and expected Lyapunov exponents associated with toss-and-catch dynamics. 
It also includes additional higher-point examples, further comparisons between the logistic IFS and logistic-tent IFS, and the algorithm used to approximate the attractor “size”.

\section*{Acknowledgements}
This work was partly supported by JSPS KAKENHI (Grant Nos. 24K15117 to T. Onozaki and 24K00537 to Y. Saiki) and JST (Grant No. JPMJSP2174 to Y. Sugita).
H. Kato gratefully acknowledges support from the Heiwa Nakajima Foundation through a scholarship.
The authors would like to thank Natsuki Tsutsumi, Miki U. Kobayashi, Yuzuru Sato, and Yoshiyuki Y. Yamaguchi for their insightful comments and for informing us of related works.


\begin{thebibliography}{14}
\bibitem [{NOT()}]{NOTE}%
  \BibitemOpen
  \href@noop {} {}\bibinfo {note} {Technically, we call this an $n$-point toss-and-catch only if every point is reachable from every initial point.}\BibitemShut {Stop}%
\end{thebibliography}
%

\end{document}


\section{Supplementary calculations for stationary measures and Lyapunov exponents}

This supplementary material provides detailed calculations of transition matrices, stationary distributions, and expected Lyapunov exponents for higher-point toss-and-catch structures and the logistic-tent IFS cases.
The toss-and-catch dynamics can be represented as finite-state Markov chains.
This allows us to determine invariant measures and calculate the expected Lyapunov exponents, denoted by $E[\lambda_p]$.

Let $P$ be the row-stochastic transition matrix of the Markov chain, where $P_{ij}$ represents the transition probability from the $i$-th point to the $j$-th point. 
We write stationary distributions as row vectors. 
The stationary distribution vector $\pi$ satisfies
\[
\pi = \pi P.
\]
By using the stationary distribution and the probability $p$ of map selection, we obtain the expected Lyapunov exponent on the toss-and-catch.
For a two-map IFS, the expected Lyapunov exponent can be written as
\[
E[\lambda_p]
=
\sum_i \pi_i
\left(
p \ln |g'(x_i)| + (1-p)\ln |h'(x_i)|
\right),
\]
provided that the relevant derivatives exist at the points of $\Lambda$.

\subsection{Logistic IFS}


\subsubsection{3-Point Case}
The set of points is:
\[
    \left\{\frac{1}{\alpha}, \frac{\alpha-1}{\alpha},  \frac{\beta+1 + \sqrt{(\beta+1)(\beta-3)}}{2\beta}\right\}.
\]
The transition matrix $P$ is:
\[
    P=
    \begin{pmatrix}
        0 & 1 - p & p\\
        0 & 1 - p & p\\
        1 - p & p & 0
    \end{pmatrix}.
\]
The stationary distribution $\pi$ is:
\[
    \pi = \frac{1}{p+1}\Big( p(1-p),\; 1 - p(1-p),\; p \Big).
\]
The expected Lyapunov exponent is approximately:
\[
    E[\lambda_p] \approx  \frac{-2.58 + p \left(-3.62 + p \left( 1.27p - 4.12 \right) \right)}{6(1+p)}.
\]
Specifically, $E[\lambda_p] \approx -0.584959$ if $p=0.5$. Numerically, we observe that $E[\lambda_p] < 0$ for all $p \in (0,1)$.

\subsubsection{5-Point Case}
The set of points is:
\begin{equation}
    \begin{split}
        \bigg\{&\frac{1}{\beta}, \frac{\beta+1 - \sqrt{(\beta+1)(\beta-3)}}{2\beta}, \frac{\alpha - 1}{\alpha}, \\
        &\frac{\beta - 1}{\beta}, \frac{\beta+1 + \sqrt{(\beta+1)(\beta-3)}}{2\beta} \bigg\}.
    \end{split}
\end{equation}

The transition matrix $P$ is:
\begin{equation}
    P=
    \begin{pmatrix}
        0 & 1 - p & 0 & p & 0\\
        0 & 0 & 1 - p & 0 & p\\
        0 & 0 & 1 - p & 0 & p\\
        0 & 1 - p & 0 & p & 0\\
        1 - p & p & 0 & 0 & 0
    \end{pmatrix}.
\end{equation}
The stationary distribution $\pi$ is:
\[
    \pi = \frac{1}{2p+1}\Big( p(1-p),\; p,\; 1-p,\; p^2,\; p \Big).
\]
The expected Lyapunov exponent is approximately:
\begin{equation}
    \begin{split}
        E[\lambda_p] \approx \frac{1.25 \times 10^{-3} p}{1 + 2p} & \Big[ \left( (12.7 + 109p) + 11.7 \right) \\
        & \times \left( 7.15p -20.3\right) \Big].
    \end{split}
\end{equation}
At $p=0.5$, $E[\lambda_p] \approx -0.473$. Numerically, we observe that $E[\lambda_p] < 0$ for all $p \in (0,1)$.
\subsection{Logistic-tent IFS}
In the following logistic-tent calculations, we consider the nontrivial case $\mu>1$.
\subsubsection{1-Point Case}
For the one-point case of logistic-tent IFS, the invariant measure is concentrated at a single fixed point. Let $p$ be the probability of selecting the logistic map with derivative $|\gamma(1-2x)|$, and $(1-p)$ be the probability of selecting the tent map with slope $|\mu|$.
The expected Lyapunov exponent is:
\begin{align*}
    E[\lambda_p]
    &= p \ln (\mu -1) + (1-p) \ln \mu.
\end{align*}
$E[\lambda_p] > 0$ when $\mu>\frac12(1+\sqrt{5})$ with $p=0.5.$
\subsubsection{2-Point Case}
For the two-point case of logistic-tent IFS, the structure mirrors the logistic IFS case, so the stationary distribution $\pi$ is:
\[
    \pi = (p, 1-p)
\]
The expected Lyapunov exponent is:
\begin{align*}
    E[\lambda_p]
    &= p \ln \frac{\mu -1}{\mu} + (1-p) \ln \mu.
\end{align*}
$E[\lambda_p] > 0$ when $\mu>2$ with $p=0.5$.

\subsubsection{3-Point Case}
The 3-point toss-and-catch set consists of the nontrivial fixed point of $g$ and the two points on a period-2 orbit of $h$.
The set of points is:
\[
    \left\{\frac{\mu}{1+\mu^2}, 1-\frac{1}{\gamma}, \frac{\mu^2}{1+\mu^2}\right\}.
\]
The transition matrix $P$ is:
\[
    P =
    \begin{pmatrix}
        0 & 0 & 1\\
        1 - p & p & 0\\
        1 - p & p & 0
    \end{pmatrix}.
\]
The stationary distribution $\pi$ is:
\[
    \pi = \frac{1}{2-p}\Big(1-p,\; p,\; 1-p\Big).
\]
The expected Lyapunov exponent, using $\gamma = \frac{1+\mu^2}{\mu}$, is:
\begin{align*}
    E[\lambda_p] &= \frac{p(3-p)}{2-p}\ln(\mu-1) \\
    &\quad + \frac{p(1-p)}{2-p}\ln(\mu+1) + (1-2p)\ln\mu.
\end{align*}
Numerically, $E[\lambda_p] > 0$ when $\mu\gtrsim 1.81$ with $p=0.5$.

\section{3-point toss-and-catch of logistic IFS and logistic-tent IFS}
The 3-point toss-and-catch of the logistic IFS possesses a bridging point, whereas the logistic-tent IFS does not. This distinction arises because the logistic and tent maps intersect at a point in the interval $(0,1)$ in addition to the trivial intersections at $x=0$ and $x=1$.
We show that logistic IFS cannot realize a 3-point toss-and-catch with a fixed point and a period-2 orbit without bridging points.

Assume $g$ and $h$ are logistic maps with different parameters. Let $x_g$ denote the nontrivial fixed point of $g$, and let $\{x_1, x_2\}$ be the period-2 orbit of $h$. Suppose that a 3-point toss-and-catch is realized on the set $S = \{x_g, x_1, x_2\}$. The set $S$ must be forward invariant under the IFS; in particular, $\{g(x_1), g(x_2)\} \subset S$. Since $x_1$ and $x_2$ constitute a period-2 orbit of $h$, they are not fixed by $h$, and by hypothesis, they are distinct from the fixed point $x_g$ of $g$.

We argue that $g$ must map both period-2 points to the fixed point $x_g$. First, $g(x_1) \neq x_1$ because $x_1\neq x_g$.
Second, $g(x_1)\neq x_2$. Indeed, since $h(x_1)=x_2$, the equality $g(x_1)=x_2$ would imply $g(x_1)=h(x_1)$. For two logistic maps with distinct parameters, this can occur only at the trivial points $x=0$ or $x=1$, whereas the period-2 points considered here lie in $(0,1)$.
The same argument applied to $x_2$ gives $g(x_2)=x_g$.
Since $x_g$ is a fixed point of $g$, it follows that
\[
g(x_g)=g(x_1)=g(x_2)=x_g.
\]

This implies that $x_g$ has at least three distinct preimages under $g$: $x_g, x_1,$ and $x_2$. This is impossible for the logistic map, which is a unimodal quadratic map and therefore admits at most two preimages for any value in its range. This contradiction demonstrates that such a 3-point toss-and-catch cannot exist for the logistic IFS, provided the points are distinct.

In contrast, the logistic-tent system avoids this obstruction. The additional intersection point of the logistic and tent maps within $(0,1)$ permits toss-and-catch transitions without forcing both $g(x_1)$ and $g(x_2)$ to collapse onto $x_g$. Specifically, it becomes possible to satisfy $g(x_1) = x_2$ consistent with the relations $h(x_1) = x_2$ and $h(x_2) = x_1$. Thus, $x_g$ need not have three distinct preimages. This mechanism allows the logistic-tent system to avoid the bridging-point phenomenon observed in the purely logistic IFS.

\section{Approximation of attractor ``size''}
To detect the toss-and-catch structure in parameter space, we estimate the size of the attractor for each parameter value by computing the minimum number of $\varepsilon$-intervals required to cover a trajectory. 
Since the systems considered in this study are composed of one-dimensional maps,
this quantity can be computed efficiently.

First, we generate a sufficiently long trajectory and sort the points in ascending order. We then apply a greedy covering procedure: starting from the leftmost uncovered point, we place an interval so that this point becomes the left endpoint of the interval. The interval has length $2\varepsilon$, which corresponds to an $\varepsilon$-interval in one dimension. We repeat this procedure until all points of the trajectory are covered (\textbf{Fig.~S1}).

The total number of intervals obtained by this procedure gives the minimal number of $\varepsilon$-intervals required to cover the trajectory in one dimension.

\begin{figure}
    \centering
    \begin{tikzpicture}[x=1.1cm,y=1cm, line cap=round, line join=round]
        \def\eps{0.85} 
        \coordinate (p1) at (1.20,0);
        \coordinate (p2) at (2.65,0);
        \coordinate (p3) at (3.55,0);
        \draw[very thick] (0,0) -- (6.8,0);
        \draw[very thick] (0,0.12) -- (0,-0.12);
        \draw[very thick] (6.8,0.12) -- (6.8,-0.12);
        \foreach \P in {p1,p2,p3}{
            \fill (\P) circle (2.8pt);
        }
        \draw[ultra thick, red!80!black]
        ($(p1)+(0,0)$) -- ($(p1)+(2*\eps,0)$);
        \draw[ultra thick, red!80!black]
        ($(p1)+(0,0.18)$) -- ($(p1)+(0,-0.18)$);
        \draw[ultra thick, red!80!black]
        ($(p1)+(2*\eps,0.18)$) -- ($(p1)+(2*\eps,-0.18)$);
        \draw[ultra thick, red!80!black]
        ($(p3)+(0,0)$) -- ($(p3)+(2*\eps,0)$);
        \draw[ultra thick, red!80!black]
        ($(p3)+(0,0.18)$) -- ($(p3)+(0,-0.18)$);
        \draw[ultra thick, red!80!black]
        ($(p3)+(2*\eps,0.18)$) -- ($(p3)+(2*\eps,-0.18)$);
        \draw[<->, thick, red!80!black]
        ($(p1)+(0,0.55)$) -- node[above, inner sep=1pt] {\allred $2\varepsilon$} ($(p1)+(2*\eps,0.55)$);
        \draw[<->, thick, red!80!black]
        ($(p3)+(0,0.55)$) -- node[above, inner sep=1pt] {\allred $2\varepsilon$} ($(p3)+(2*\eps,0.55)$);
        \node[below=4pt] at (0, 0) {$0$};
        \node[below=4pt] at (p1) {$x_{(1)}$};
        \node[below=4pt] at (p2) {$x_{(2)}$};
        \node[below=4pt] at (p3) {$x_{(k)}$};
        \node[below=4pt] at (6.8, 0) {$1$};
    \end{tikzpicture}
    \caption{\textbf{Greedy $\varepsilon$-cover in one dimension.}
        Given a sorted trajectory $\{x_{(i)}\}$, start from the leftmost uncovered point and place an interval
        $[x_{(i)},x_{(i)}+2\varepsilon]$ (equivalently an $\varepsilon$-interval centered at $x_{(i)}+\varepsilon$). Repeat until all points are covered. The number of intervals is the estimated cover number. [ALT Text: Schematic illustration of the greedy $\varepsilon$-cover algorithm used to estimate attractor size in one dimension. Sorted trajectory points are sequentially covered by intervals of length 2$\varepsilon$, and the total number of intervals gives the estimated cover number.]}
    \label{fig:eps-cover-greedy}
\end{figure}
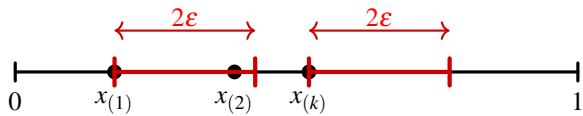

\bibliography{bib/functionmap, bib/KOSS-IFS}